\documentclass[11pt]{article}
\usepackage{mathrsfs}
\usepackage{amssymb}
\usepackage{amsmath}
\usepackage{amsbsy}
\usepackage{bbold}
\usepackage{epsfig}
\usepackage{enumerate}
\usepackage{bm}
\usepackage{color}
\usepackage{psfrag}
\usepackage{type1cm}
\usepackage[parfill]{parskip}
\usepackage{float}
\usepackage{graphicx}
\usepackage{epsfig}
\usepackage{amsmath, amsthm, amssymb}
\usepackage{verbatim}
\usepackage{multicol}
\usepackage{mathptmx} 
\usepackage{courier}
\usepackage{url}
\usepackage{latexsym}
\usepackage{mathrsfs}
\usepackage[colorlinks = true,
            linkcolor  = blue,
            citecolor  =blue,
            urlcolor   =blue,
            bookmarks  = true]{hyperref}
 \usepackage{enumitem}
\setlist[enumerate,1]{label=(\roman*)}
\usepackage[normalem]{ulem}
\usepackage{bm}
\usepackage{dsfont}
\usepackage{stmaryrd}
\usepackage[dvipsnames]{xcolor}

\usepackage[backend=biber, style=numeric, citestyle=numeric-comp]{biblatex}
\DeclareFieldFormat{labelnumberwidth}{\mkbibbrackets{#1}}
\DeclareNameAlias{author}{family-given}

\newtheorem{definition}{ \noindent D{\footnotesize EFINITION}}[section]
\newtheorem{theorem}{ \noindent T{\footnotesize HEOREM}}

\newtheorem{lemma}{ \noindent L{\footnotesize EMMA}}[section]
\newtheorem{coro}{ \noindent C{\footnotesize OROLLARY}}
\newtheorem{remark}{ \noindent R{\footnotesize EMARK}}[section]

\newcommand{\E}{\mathbb{E}}
\newcommand{\Prob}{\mathbb{P}}

\DeclareFieldFormat*{title}{#1}

\AtEveryBibitem{%
  \clearfield{issn}
  \clearfield{url}
  \clearfield{eprint}
}

\renewbibmacro{in:}{}

\DeclareFieldFormat[article]{volume}{\textbf{#1}}
\DeclareFieldFormat[article]{number}{\textbf{#1}}
\renewbibmacro*{volume+number+eid}{%
  \printfield{volume}
  \iffieldundef{number}
    {}%
    {\printtext[parens]{\printfield{number}}}
  \setunit{\addcomma\space}%
  \printfield{eid}}

\DeclareFieldFormat[journal]{journaltitle}{\mkbibemph{#1}}

\DeclareNameAlias{author}{family-given}
\DeclareNameAlias{editor}{family-given}

\DeclareDelimFormat{namedelim}{\addcomma\space}
\DeclareFieldFormat[article,incollection,inproceedings,book]{giveninit}{\mkbibnamegiven{#1}}

\DeclareNameFormat{family-given}{%
  \usebibmacro{name:family-given}
    {\namepartfamily}
    {\namepartgiveni} 
    {\namepartprefix}
    {\namepartsuffix}%
  \usebibmacro{name:andothers}}

\DeclareFieldFormat{pages}{#1}

\numberwithin{equation}{section}

\theoremstyle{plain}

\def\author#1{\par
    {\centering{\authorfont#1}\par\vspace*{0.05in}}
}

\def\titlefont{\fontsize{13}{15}\bfseries\boldmath\selectfont\centering{}}
\def\authorfont{\fontsize{13}{15}}

\def\title#1{
    \thispagestyle{plain}
    \vspace*{-14pt}
    \vskip 79pt
    {\centering{\titlefont #1\par}}%
    \vskip 1em
}

\renewcommand{\baselinestretch}{1.5}

\begin{document}

\begin{center}{\LARGE\bf The exact group-sparse recovery for block diagonal matrices with subexponential entries}
\end{center}

\begin{center} {\sc
Guozheng Dai\footnote{Department of Mathematics, Hong Kong University of Science and Technology;}, Tiankun Diao\footnote{Corresponding author, Institute for Financial Studies, Shandong University,  tiankundiao@gmail.com},  Hanchao Wang\footnote{Institute for Financial Studies, Shandong University}.}
\end{center}

 \renewcommand{\abstractname}{~}
\begin{abstract}
 {\bf Abstract:}   
 We study block-diagonal random matrices with i.i.d. subexponential entries and show that, despite their highly structured form, they already guarantee exact sparse recovery from a nearly optimal number of measurements. When the matrix reduces to a single block, our framework collapses to the classical i.i.d. subexponential ensemble, and our bounds recover the well-known optimal rates previously established for unstructured random matrices.



\vspace{-3mm}
\end{abstract}

\renewcommand{\baselinestretch}{1.2}


\section{Introduction}\label{sec:intro}
A common problem in modern signal processing applications is the recovery of sparse signals from incomplete sampling data. Specifically, the goal is to reconstruct an unknown $s$-sparse signal $x \in \mathbb{C}^n$, where the $s$-sparsity is defined as $\|x\|_0 := |\{l : x_l \neq 0\}| \leq s$, from the observed linear measurements $y = \Phi x \in \mathbb{C}^m$. Here, $\Phi$ is a known $m \times n$ matrix, referred to as the measurement matrix. The problem becomes particularly interesting when $m \ll n$, rendering the system underdetermined. In this scenario, reconstructing $x$ would typically be impossible without additional information. However, if it is known a priori that $x$ is $s$-sparse, the situation changes significantly.

A widely applied approach for reconstruction is $\ell_{1}$-minimization \cite{Candes_ieee_2006}
\begin{align}\label{Eq_intro_l1_minimization}
    \text{min}_{z\in\mathbb{C}^{n}} \Vert z\Vert_{1}\quad \quad \text{subject to}\,\, \Phi z=y,
\end{align}
where $\Vert z\Vert_{p}=(\sum_{i}\vert z_{i}\vert^{p})^{1/p}$ denotes the usual $\ell_{p}$-norm.  A series of landmark papers have demonstrated that any $s$-sparse vector $x \in \mathbb{C}^n$ can be exactly reconstructed from the linear measurements $y = \Phi x$, where $\Phi \in \mathbb{C}^{m \times n}$, and $m \ll n$, using $\ell_1$-minimization. This reconstruction is possible as long as the measurement matrix $\Phi$ satisfies certain geometric properties. For example, Cand\'{e}s and Tao \cite{Candes_ieee_2005} (see also \cite{Mendelson_gafa,Vershynin_high_dimension}) proved that if the measurement matrix \(\Phi\) is a random matrix with i.i.d. (independent and identically distributed) subgaussian entries (see Definition \ref{Def_subgaussian} below) having mean \(0\) and variance \(1\), then under the condition
\begin{align}\label{Eq_intro_condition_rip}
    m\ge Cs\log (en/s),
\end{align}
one can exactly reconstruct the sparse signal $x$ via $\ell_{1}$-minimization with high probability. Here $C>0$ is an absolute constant and $e$ is the exponential constant. Note that the condition \eqref{Eq_intro_condition_rip} is optimal in the sense that if it is not satisfied, the recovery fails with high probability \cite[Chapter 10]{Foucart_book}.  Moreover, analogous results have been obtained for heavier-tailed ensembles, such as measurement matrices whose entries are independent copies of subexponential random variables \cite{Adamczak_Constr,Foucart_LAA,Foucart_Studiamath}. 

While random i.i.d. matrices are highly desirable from a theoretical standpoint, practical applications often demand additional structure in the measurement matrix \(\Phi\). On the one hand, in most engineering applications, the structural properties of the measurement system are largely predetermined by the specific requirements of the application itself. On the other hand, fast matrix multiplication is not available for such unstructured matrices, which makes the recovery algorithms run very slowly. Typical examples of structured random matrices are partial Fourier transform matrices \cite{Rudelson_cpam}, partial circulant matrices \cite{Krahmer_cpam}, and so on.

As observed above, all the examples of measurement matrices we discussed are random matrices. It is important to note that the deterministic construction of measurement matrices with provably optimal scaling in terms of the information dimension of signals (see condition \eqref{Eq_intro_condition_rip}) remains an unsolved problem.

In this paper, we consider a special class of structured random matrices, block diagonal measurement matrices whose blocks are independent random matrices. In particular, let $\Phi_{1},\cdots, \Phi_{L}$ be a sequence of independent $m\times d$ random matrices with i.i.d. entries. We are interested in the following $mL\times dL$ structured random matrix:
\begin{align}\label{Eq_intro_Phi}
    \Phi = \begin{pmatrix} \Phi_1 & & \\ & \ddots & \\ & & \Phi_L \end{pmatrix}.
\end{align}
Such block-wise measurement models appear in various applications, such as distributed compressed sensing \cite{Duarte_conference} and the multiple measurement vector model \cite{Chenjie_ieee_signal}. Recently, Koep et al. \cite{Koep_acha} demonstrated that when the nonzero entries of \(\Phi\) are subgaussian random variables, it is possible to exactly reconstruct a class of more structured signals, specifically those whose nonzero coefficients appear in groups. Next, we will give a specific expression of the results in \cite{Koep_acha}.
Before we proceed, we recommend that readers refer to Section \ref{sec:notation} to familiarize themselves with the notations used there.

\begin{theorem}[Theorem 3.2 in \cite{Koep_acha}]\label{Theo_1}
    Let  $\mathcal{I}=\{\mathcal{I}_{1}, \cdots, \mathcal{I}_{G}\}$ be a group partition of $[dL]$, and let $\Psi$ be a $dL\times dL$ unitary matrix. Denote $g:=\max_{i\in[G]}\vert \mathcal{I}_{i}\vert$. Let $\Phi$ be an $mL\times dL$ random matrix defined in \eqref{Eq_intro_Phi} and assume that the nonzero entries of $\Phi$ are i.i.d. subgaussian random variables with mean $0$, variance $1$, and subgaussian norm $\tau$. Then for the following linear program
    \begin{align}
        y=\Phi\Psi x, \quad x\in \Sigma_{\mathcal{I}, s}(\mathbb{C}^{dL}),
    \end{align}
if 
\begin{align}\label{Eq_condition_theorem1}
    m\ge C(\tau)\max\Big( s \mu_{\mathcal{I}}^{2}\log(dL)\log^{2}s\big(\log G+g\log(s/\mu_{\mathcal{I}}) \big), \log \eta^{-1}\Big),
\end{align}
where $\mu_{\mathcal{I}}=\min\Big\{ \sqrt{d}\max_{i\in [dL]}\Vert \Psi_{i,:}^{\top}\Vert_{\mathcal{I}, \infty}, 1 \Big\}$, with probability at least $1-\eta$, one can exactly recover the $s$-group sparse vector $x$ via the following group $\ell_{1}$-minimization:
\begin{align}
    \min_{z\in\mathbb{C}^{n}} \Vert z\Vert_{\mathcal{I}, 1}\quad \quad \textnormal{subject to}\,\,  Az=y.
\end{align}
\end{theorem}

Our main result concentrates on the block diagonal matrices with subexponential entries. We show that such structured random matrices can lead to the exact sparse recovery with the nearly optimal number of measurements.

\begin{theorem}\label{Theo_main}
    Let  $\mathcal{I}=\{\mathcal{I}_{1}, \cdots, \mathcal{I}_{G}\}$ be a group partition of $[dL]$, and let $\Psi$ be a $dL\times dL$ unitary matrix. Denote $g:=\max_{i\in[G]}\vert \mathcal{I}_{i}\vert$. Let $\Phi$ be an $mL\times dL$ random matrix defined in \eqref{Eq_intro_Phi} and assume that the nonzero entries of $\Phi$ are i.i.d. subexponential random variables with mean $0$, variance $1$, and subexponential norm $\tau$. Then for the following linear program
    \begin{align}
        y=\Phi\Psi x, \quad x\in \Sigma_{\mathcal{I}, s}(\mathbb{C}^{dL}),
    \end{align}
if 
\begin{align}\label{Eq_condition_theorem2}
        m\ge C(\tau)\max\Big( gLs\log\big(\frac{dL}{gs} \big), \log \eta^{-1}\Big),
    \end{align}
 with probability at least $1-\eta$, one can exactly recover the $s$-group sparse vector $x$ via the following group $\ell_{1}$-minimization:
\begin{align}
    \min_{z\in\mathbb{C}^{n}} \Vert z\Vert_{\mathcal{I}, 1}\quad \quad \textnormal{subject to}\,\,  Az=y.
\end{align}
\end{theorem}

\begin{remark}
    When $L=1$, $\mathcal{I}=\{\{1\},\cdots,\{d\}  \}$, and $\Psi$ is an identity matrix, i.e. $\Phi\Psi$ is an $m\times d$ random matrix with i.i.d. subexponential entries, the condition \eqref{Eq_condition_theorem2} changes into
    \begin{align}
        m\ge C(\tau)\max\Big( s\log\big(\frac{d}{s} \big), \log \eta^{-1}\Big),
    \end{align}
    which is the optimal condition for exact recovery via $\ell_{1}$-minimization.
    Hence, Theorem \ref{Theo_main} recovers the classical result, see \cite[Theorem 6.1]{Foucart_Studiamath} for details.
\end{remark}
\begin{remark}
    As noted in \cite[Section 7]{Koep_acha}, there is an additional \(\log^{2}s\) factor in \eqref{Eq_condition_theorem1}. The authors of \cite{Koep_acha} conjectured that this factor should be eliminated. Theorem \ref{Theo_main} confirms that this is indeed feasible.
\end{remark}

\section{Preliminaries}

\subsection{Notations}\label{sec:notation}

Unless otherwise specified, we use $C, C_1, C_2, c, c_1, c_2, \ldots$ to denote universal constants that are independent of matrix dimensions and the parameters of random variables. Similarly, $C(\varepsilon), C_1(\varepsilon), \ldots$ represent constants that depend only on the parameter $\varepsilon$. These constants may change from line to line.

For an integer $n\in \mathbb{N}$, we use the notation $[n]:=\{1, \cdots, n\}$. Let $\mathcal{I}:=\{\mathcal{I}_{1}, \cdots, \mathcal{I}_{m}  \}$ be a group partition of $[n]$, namely $\mathcal{I}_{i}\cap\mathcal{I}_{j}=\emptyset$ for all $1\le i< j\le m$, and $\cup_{i=1}^{m}\mathcal{I}_{i}=[n]$. 

Given a vector \( a = (a_1, \cdots, a_n)^\top \in \mathbb{C}^n \), for \( p \ge 1 \), the \( \ell_p \)-norm of \( a \) is defined as \( \|a\|_p = \left( \sum_{i=1}^n |a_i|^p \right)^{1/p} \). We also denote $\Vert a\Vert_{0}=\vert \{i: \vert a_{i}\vert\neq 0  \}\vert$. For a subset \( I \subset [n] \), the vector \( a_I \in \mathbb{C}^n \) is such that the coordinates indexed by \( I \) are the same as those in \( a \), and all other coordinates are zero. Given a group partition $\mathcal{I}:=\{\mathcal{I}_{1}, \cdots, \mathcal{I}_{m} \}$, for $p\ge 1$, the $\ell_{\mathcal{I}, p}$-norm of $a$ is defined as follows:
\begin{align}
    \Vert a\Vert_{\mathcal{I}, p}:=\Big(\sum_{i=1}^{m} \Vert a_{\mathcal{I}_{i}}\Vert_{2}^{p} \Big)^{1/p}.
\end{align}
The $\ell_{\mathcal{I}, 0}$-norm of $a$ is defined as follows:
\begin{align}
    \Vert a\Vert_{\mathcal{I}, 0}:=\big\vert\{i\in [m]:a_{\mathcal{I}_{i}}\neq (0,0,\cdots,0)^\top \} \big \vert.
\end{align}
  We usually use the following $s$-group-sparse vector sets:
\begin{align}
    \Sigma_{\mathcal{I}, s}(\mathbb{C}^{n})=\Big\{a\in\mathbb{C}^{n}: \Vert a\Vert_{\mathcal{I}, 0}\le s  \Big\}.
\end{align}

Given an \( m \times n \) matrix \( A \), we denote the \( i \)-th column of \( A \) by \( A_{:,i} \) and the \( i \)-th row of \( A \) by \( A_{i,:} \). Let $I\subset [n]$, then $A_{:,I}=(A_{:, i}, i\in I)$ is an $m\times \vert I\vert$ submatrix of $A$. Similarly, given a set $I\subset [m]$, $A_{I, :}=(A_{i, :}^{\top}, i\in I)^{\top}$ is an $\vert I\vert\times n$ submatrix of $A$.

\subsection{$\alpha$-subexponential random variables}

For $0<\alpha\le2$, the concept of $\alpha$-subexponential random variable is as follows:
\begin{definition}\label{Def_subgaussian}
    A random variable $\xi$ is called $\alpha$-subexponential if its tail probability satisfies for $t\ge 0$
    \begin{align}
        \Prob\{\vert \xi\vert\ge t  \}\le 2\exp\big(-(\frac{t}{K})^{\alpha}\big),
    \end{align}
    where $K>0$ is a parameter. The variable \(\xi\) is typically referred to as a subgaussian variable when \(\alpha = 2\) and as a subexponential variable when \(\alpha = 1\). The \(\alpha\)-subexponential norm of \(\xi\) is defined as follows:
    \begin{align}
        \Vert \xi\Vert_{\Psi_{\alpha}}:=\inf\big\{t>0: \E\exp\big(\frac{\vert \xi\vert^{\alpha}}{t^{\alpha}}  \big)\le 2  \big\}.
    \end{align}
\end{definition}

Then, we present a sufficient condition for a subexponential random variable.
\begin{lemma}[Proposition 2.7.1 in \cite{Vershynin_high_dimension}]\label{Lem_subexponential_p_moment}
    Let $\xi$ be a  random variable such that, for $p\ge 1$
    \begin{align}
        (\mathbb{E}\vert\xi\vert^{p})^{1/p}\le Kp.
    \end{align}
    Then, $\xi$ is a subexponential random variable with $K/C\le \Vert \xi\Vert_{\Psi_{1}}\le CK$, where $C>0$ is an absolute constant.
\end{lemma}

Next, we introduce the concentration property of sums of independent subexponential random variables, which is known as Bernstein's inequality. For the general case of $\alpha$-subexponential random variables, we refer interested readers to \cite{Gluskin_Studiamath,Hitczenko_Studiamath} and do not cover it here.

\begin{lemma}[Theorem 2.8.1 in \cite{Vershynin_high_dimension}]\label{Lem_berstein}
    Let $\xi_{1}, \cdots, \xi_{n}$ be independent mean-zero sub-exponential random variables. Then, for every $t\ge 0$, we have
    \begin{align}
        \Prob \Big\{ \big\vert\sum_{i=1}^{n}\xi_{i} \big\vert \ge t \Big\}\le 2\exp\Big(-c\min\big(\frac{t^{2}}{\sum_{i=1}^{n}\Vert \xi_{i}\Vert_{\Psi_{1}}^{2}}, \frac{t}{\max_{i}\Vert \xi_{i}\Vert_{\Psi_{1}}}  \big)  \Big),
    \end{align}
    where $c>0$ is an absolute constant.
\end{lemma}

\subsection{$\varepsilon$-net}
Consider a normed space $(X, \Vert \cdot\Vert)$.
We call a subset $\mathcal{N}\subset U$ an $\varepsilon$-net of $U$ if it satisfies that
\begin{align}
	\inf_{x\in \mathcal{N}}\Vert x-u\Vert <\varepsilon, \qquad \forall u\in U.\nonumber
\end{align}

\begin{lemma}[Corollary 4.2.13 in \cite{Vershynin_high_dimension}]\label{Lem_epsilon_net}
Let $\varepsilon\in (0, 1)$ and consider the unit Euclidean ball $B_{2}^{n}$ in $\mathbb{R}^{n}$, i.e. $B_{2}^{n}=\{x\in\mathbb{R}^{n}: \Vert x\Vert_{2}\le 1\}$. Then there exists an  $\varepsilon$-net $\mathcal{N}$ such that 
\begin{align}
    \vert \mathcal{N}\vert \le \big(\frac{3}{\varepsilon}\big)^{n}.
\end{align}
\end{lemma}
 Based on Lemma \ref{Lem_epsilon_net}, one can immediately bound the cardinality of the $\varepsilon$-net of  the unit Euclidean ball in $\mathbb{C}^{n}$, reading as follows:
\begin{coro}\label{Coro_cardinality_unitball}
    For the unit Euclidean ball in $\mathbb{C}^{n}$ and $\varepsilon\in (0, 1)$, there exists an $\varepsilon$-net $\mathcal{N}$ such that 
    \begin{align}
    \vert \mathcal{N}\vert \le \big(\frac{3}{\varepsilon}\big)^{2n}.
\end{align}
\end{coro}

\begin{coro}\label{Coro_epsilon_net}
Let $\Omega := \{x \in \mathbb{C}^n : \|x\|_2 = 1, \|x\|_0 \leq s\}$. Then there exists an $\varepsilon$-net  $\mathcal{N}$ of $\Omega$ such that
$$|\mathcal{N}| \leq\left(\frac{9e  n}{\varepsilon^{2}  s}\right)^s.$$
\end{coro}
\begin{proof}
Denote $\Omega_{1}:=\{x\in\mathbb{C}^{n}: x_{i}=0, 1\le i\le n-s  \}$. Hence, Corollary \ref{Coro_cardinality_unitball} implies that there exists an $\varepsilon$-net $\mathcal{N}_{1}$ of $\Omega\cap\Omega_{1}$ such that
\begin{align}
    \vert \mathcal{N}_{1}\vert\le \left(\frac{3}{\varepsilon}\right)^{2s}.\nonumber
\end{align}

Thus, we can construct an $\varepsilon$-net $\mathcal{N}$ for $\Omega$ with cardinality satisfying
$$
|\mathcal{N}| \leq \binom{n}{s} \left(\frac{3}{\varepsilon}\right)^{2s}\le  \left(\frac{e  n}{s}\right)^{s}  \left(\frac{3}{\varepsilon}\right)^{2s} = \left(\frac{9e  n}{\varepsilon^{2}  s}\right)^s.
$$   
\end{proof}

\section{Restricted isometry property}
In this section, we first introduce a well-established tool for analyzing the performance of sparse recovery, known as the restricted isometry property (RIP). An \( m \times n \) matrix \( A \) is said to satisfy the RIP with parameters \( \alpha \), \( \beta \), and \( s \) if the inequality
\[
\alpha \| x \|_2 \leq \| A x \|_2 \leq \beta \| x \|_2
\]
holds for all \( s \)-sparse vectors \( x \) of dimension \( n \). Cand\'{e}s and Tao \cite{Candes_ieee_2005} proved that if the measurement matrix satisfies RIP, then one can recover the sparse vector via \(\ell_{1}\)-minimization. In fact, many works (see \cite{Candes_ieee_2005,Vershynin_high_dimension}) have demonstrated that a random i.i.d. subgaussian matrix can be used as a measurement matrix for exact sparse recovery by virtue of proving that it has the RIP.

However, for random i.i.d. subexponential matrices, using RIP to analyze their performance in sparse recovery does not seem to be appropriate. This is because these matrices do not satisfy RIP with the optimal number of measurements (see the condition \eqref{Eq_intro_condition_rip}). Foucart and Lai \cite{Foucart_Studiamath} introduced an alternative tool for analyzing sparse recovery, known as the restricted isometry property with respect to the \(\ell_{1}\) and \(\ell_{2}\) norms (RIP\(_{1,2}\)). An \(m \times n\) matrix \(A\) is said to satisfy the RIP\(_{1,2}\) with parameters \(\alpha\), \(\beta\), and \(s\) if the inequality
\[
\alpha \| x \|_2 \leq \| A x \|_1 \leq \beta \| x \|_2
\]
holds for all \(s\)-sparse vectors \(x\) of dimension \(n\). Foucart and Lai \cite{Foucart_Studiamath} proved that the RIP\(_{1,2}\) is a sufficient condition for exact sparse recovery. They also showed that random i.i.d. subexponential matrices satisfy the RIP\(_{1,2}\) with the optimal number of measurements.

To address the problems of group sparse recovery, we first introduce the following tool,  a modified version of RIP$_{1,2}$.

 \begin{definition}
    Given a group partition \(\mathcal{I}\) of \([n]\), an \(m \times n\) matrix \(A\) is said to satisfy the group restricted isometry property with respect to the \(\ell_{1}\) and \(\ell_{2}\) norms (group-RIP$_{1, 2}$) with parameters \(\alpha\), \(\beta\), and \(s\) if
\[
\alpha \| x \|_{2} \leq \| A x \|_{1} \leq \beta \| x \|_{2}, \quad \forall x \in \Sigma_{\mathcal{I}, s}(\mathbb{C}^{n}).
\]
 \end{definition}
Next, we will demonstrate that if the measurement matrix satisfies group-RIP$_{1, 2}$, the group sparse vector can be recovered via group $\ell_{1}$-minimization.
\begin{theorem}\label{Theo_rip_12}
    Suppose that an $m\times n$ matrix $A$ fulfills the group-RIP$_{1, 2}$ with certain parameters $\alpha, \beta$, and $(1+\lambda)s$, where $\lambda>(\beta/\alpha)^{2}$, namely
\[
\alpha \| x \|_{2} \leq \| A x \|_{1} \leq \beta \| x \|_{2}, \quad \forall x \in \Sigma_{\mathcal{I}, (1+\lambda)s}(\mathbb{C}^{n}).
\]
Let $\mathcal{I}:=\{ \mathcal{I}_{1}, \cdots, \mathcal{I}_{G}\}$ be a group partition of $[n]$. For the following linear program
\[
y=Ax,\quad x\in \Sigma_{\mathcal{I}, s}(\mathbb{C}^{n}),
\]
the $s$-group sparse vector $x$ can be exactly recovered via the group $\ell_{1}$-minimization:
\begin{align}\label{Eq_section3_l1_minimization}
   \min_{z\in\mathbb{C}^{n}} \Vert z\Vert_{\mathcal{I}, 1}\quad \quad \textnormal{subject to}\,\,  Az=y. 
\end{align}
\end{theorem}
\begin{proof}
    Denote by $z^{*}$ the solution of the group $\ell_{1}$-minimization \eqref{Eq_section3_l1_minimization}. Let $h=z^{*}-x$. Hence, it is enough to prove that $h=(0,0,\cdots, 0)^{\top}\in \mathbb{C}^{n}$.

We start by decomposing the set $[G]$ into several parts. Let $I_{0}\subset [G]$ be the set such that
\begin{align}
    \vert I_{0}\vert=s,\quad \text{and}\quad x_{\mathcal{I}_{j}}=(0,0,\cdots, 0)^{\top}\in \mathbb{C}^{n},\,\, \forall j\notin I_{0}.
\end{align}
Next, let $I_{1}$ be the subset of $[G]\setminus I_{0}$ such that
\begin{align}
    \vert I_{1}\vert=\lambda s,\quad \text{and}\quad \min_{i\in I_{1}}\Vert x_{\mathcal{I}_{i}}\Vert_{2}\ge \max_{i\in [G]\setminus (I_{0}\cup I_{1})}\Vert x_{\mathcal{I}_{i}}\Vert_{2}.
\end{align}
Then, denote by $I_{2}$ the subset of $[G]\setminus (I_{0}\cup I_{1})$ such that
\begin{align}
    \vert I_{2}\vert=\lambda s,\quad \text{and}\quad \min_{i\in I_{2}}\Vert x_{\mathcal{I}_{i}}\Vert_{2}\ge \max_{i\in [G]\setminus (\cup_{0\le j\le 2} I_{j})}\Vert x_{\mathcal{I}_{i}}\Vert_{2}.
\end{align}
The sets $I_{3},I_{4}, \cdots$ are defined in a similar manner.

For convenience, denote $I_{0, 1}=I_{0}\cup I_{1}$. By the definition of $h$, we have
\begin{align}\label{Eq_Theorem3.1_1}
    0=\Vert Ah\Vert_{1}=\Vert A(\sum_{i\in I_{0, 1}}h_{\mathcal{I}_{i}}+\sum_{i\in I_{0, 1}^{c}}h_{\mathcal{I}_{i}})\Vert_{1}\ge \Vert A(\sum_{i\in I_{0, 1}}h_{\mathcal{I}_{i}})\Vert_{1}-\Vert A(\sum_{i\in I_{0, 1}^{c}}h_{\mathcal{I}_{i}})\Vert_{1}.
\end{align}
Note that, by the group-RIP$_{1, 2}$ of $A$ and the fact $\vert I_{0,1}\vert\le (1+\lambda)s$, we have
\begin{align}
    \Vert A(\sum_{i\in I_{0, 1}}h_{\mathcal{I}_{i}})\Vert_{1}\ge \alpha \Vert \sum_{i\in I_{0, 1}}h_{\mathcal{I}_{i}}\Vert_{2}
\end{align}
and 
\begin{align}
    \Vert A(\sum_{i\in I_{0, 1}^{c}}h_{\mathcal{I}_{i}})\Vert_{1}\le \sum_{j\ge 2}\Vert A(\sum_{i\in I_{j}}h_{\mathcal{I}_{i}})\Vert_{1}\le  \beta\sum_{j\ge 2}\Vert \sum_{i\in I_{j}}h_{\mathcal{I}_{i}}\Vert_{2}.
\end{align}
Hence, we have by \eqref{Eq_Theorem3.1_1}
\begin{align}\label{Eq_Theorem3.1_3}
    \alpha \Vert \sum_{i\in I_{0, 1}}h_{\mathcal{I}_{i}}\Vert_{2}\le \beta\sum_{j\ge 2}\Vert \sum_{i\in I_{j}}h_{\mathcal{I}_{i}}\Vert_{2}.\nonumber
\end{align}

Note that for $j\ge 2$
\begin{align}
    \Vert\sum_{i\in I_{j}}h_{\mathcal{I}_{i}}\Vert_{2}=(\sum_{i\in I_{j}}\Vert h_{\mathcal{I}_{i}}\Vert_{2}^{2})^{1/2}&\le \Big(\big(\sum_{i\in I_{j-1}}\Vert h_{\mathcal{I}_{i}}\Vert_{2}\big)\frac{\sum_{i\in I_{j-1}}\Vert h_{\mathcal{I}_{i}}\Vert_{2}}{\lambda s}\Big)^{1/2}\\
    &\le \frac{1}{\sqrt{\lambda s}}\sum_{i\in I_{j-1}}\Vert h_{\mathcal{I}_{i}}\Vert_{2}.
\end{align}
Hence, we have
\begin{align}\label{Eq_Theorem3.1_2}
    \sum_{j\ge 2}\Vert \sum_{i\in I_{j}}h_{\mathcal{I}_{i}}\Vert_{2}\le \frac{1}{\sqrt{\lambda s}}\sum_{j\ge 1}\sum_{i\in I_{j}}\Vert h_{\mathcal{I}_{i}}\Vert_{2}.
\end{align}

Note that 
\begin{align}
    \Vert z^{*}\Vert_{\mathcal{I}, 1}=\Vert h+x\Vert_{\mathcal{I}, 1}&=\sum_{j\ge 1}\sum_{i\in I_{j}}\Vert h_{\mathcal{I}_{i}}+x_{\mathcal{I}_{i}}\Vert_{2}+\sum_{i\in I_{0}}\Vert h_{\mathcal{I}_{i}}+x_{\mathcal{I}_{i}}\Vert_{2}\\
    &=\sum_{j\ge 1}\sum_{i\in I_{j}}\Vert h_{\mathcal{I}_{i}}\Vert_{2}+\sum_{i\in I_{0}}\Vert h_{\mathcal{I}_{i}}+x_{\mathcal{I}_{i}}\Vert_{2}\\
    &\ge \sum_{j\ge 1}\sum_{i\in I_{j}}\Vert h_{\mathcal{I}_{i}}\Vert_{2}+\sum_{i\in I_{0}}\Vert x_{\mathcal{I}_{i}}\Vert_{2}-\sum_{i\in I_{0}}\Vert h_{\mathcal{I}_{i}}\Vert_{2}\\
    &=\sum_{j\ge 1}\sum_{i\in I_{j}}\Vert h_{\mathcal{I}_{i}}\Vert_{2}-\sum_{i\in I_{0}}\Vert h_{\mathcal{I}_{i}}\Vert_{2}+\Vert x\Vert_{\mathcal{I}, 1}.
\end{align}
By the definition of $z^{*}$, we have $\Vert z^{*}\Vert_{\mathcal{I}, 1}\le \Vert x\Vert_{\mathcal{I}, 1}$. Hence, we have
\begin{align}
    \sum_{j\ge 1}\sum_{i\in I_{j}}\Vert h_{\mathcal{I}_{i}}\Vert_{2}\le \sum_{i\in I_{0}}\Vert h_{\mathcal{I}_{i}}\Vert_{2}.
\end{align}
Then, \eqref{Eq_Theorem3.1_2} is further bounded by
\begin{align}
    \frac{1}{\sqrt{\lambda s}}\sum_{i\in I_{0}}\Vert h_{\mathcal{I}_{i}}\Vert_{2}\le \frac{1}{\sqrt{\lambda}}(\sum_{i\in I_{0}}\Vert h_{\mathcal{I}_{i}}\Vert_{2}^{2})^{1/2}.
\end{align}
Substituting this bound into \eqref{Eq_Theorem3.1_3}, we have
\begin{align}
    \alpha \Vert \sum_{i\in I_{0}}h_{\mathcal{I}_{i}}\Vert_{2}\le \frac{\beta}{\sqrt{\lambda}}\Vert \sum_{i\in I_{0}}h_{\mathcal{I}_{i}}\Vert_{2}.
\end{align}
Note that $\lambda>(\beta/\alpha)^{2}$, hence we have
\begin{align}
    \Vert \sum_{i\in I_{0}}h_{\mathcal{I}_{i}}\Vert_{2}=0,
\end{align}
which completes the proof.
\end{proof}

\section{Block diagonal matrices with subexponential entries}
Recall the following $mL\times dL$ block diagonal random matrix:
\begin{align}
    \Phi = \begin{pmatrix} \Phi_1 & & \\ & \ddots & \\ & & \Phi_L \end{pmatrix},
\end{align}
where $\Phi_{1}, \cdots, \Phi_{L}$ are a sequence of independent $m\times d$ random matrices with i.i.d. entries. Let $\Psi$ be a $dL\times dL$ unitary matrix. We next prove that the random matrix $\frac{1}{m}\Phi\Psi$ satisfies group-RIP$_{1,2}$ when the nonzero entries of $\Phi$ are subexponential random variables.

\begin{theorem}\label{Theo_matrix_satisfies_rip}
    Let $\mathcal{I}=\{\mathcal{I}_{1}, \cdots, \mathcal{I}_{G}\}$ be a group partition of $[dL]$ and denote $g=\max_{i\in [G]}\vert \mathcal{I}_{i}\vert$. Assume that the nonzero entries of $\Phi$ are i.i.d. subexponential random variables with mean $0$, variance $1$, and subexponential norm $\tau$. Then under the condition
    \begin{align}
        m\ge C_{1}gs\log\big(\frac{dL^{2}}{gs} \big),
    \end{align}
with the probability at least $1-e^{-cm}$, we have
\begin{align}
    C_{2}\Vert y\Vert_{2}\le \frac{1}{m}\Vert\Phi\Psi y\Vert_{1}\le C_{3}\sqrt{L}\Vert y\Vert_{2},\quad \forall y\in\Sigma_{\mathcal{I}, s}(\mathbb{C}^{dL}),
\end{align}
where $C_{1}, C_{2}, C_{3}, c$ are positive constants depending only on the subexponential norm $\tau$.
    
\end{theorem}
To facilitate reading, we shall divide the proof of Theorem \ref{Theo_matrix_satisfies_rip} into three parts, which are presented in the following three subsections in the form of lemmas.

\subsection{The concentration inequality of $\Vert\Phi\Psi x\Vert_{1}$}

\begin{lemma}\label{Lem_concentration}
 Let $\Phi, \Psi$ be the matrices defined in Theorem \ref{Theo_matrix_satisfies_rip}.   Given a fixed vector $x\in\mathbb{C}^{dL}$ satisfying $\Vert x\Vert_{2}=1$, we have for $t\ge 0$
    \begin{align}
        \mathbb{P}\left\{\Big\vert \Vert \Phi\Psi x\Vert_{1}-\mathbb{E}\Vert \Phi\Psi x\Vert_{1} \Big\vert \ge t\right\}\le 2 \exp\big(-c\min\big(\frac{t^{2}}{m}, t  \big)   \big),
    \end{align}
    where $c>0$ is a constant depending only on the subexponential norm $\tau$. 
\end{lemma}
\begin{proof}
    For convenience, let $\Psi_{1}=\Psi_{[1:d],:}, \Psi_{2}=\Psi_{[d+1 : 2d],:}, \cdots, \Psi_{L}=\Psi_{[d(L-1)+1: dL], :}$. Then, we have
    \begin{align}\label{Eq_Lemma4.1_1}
        &\mathbb{P}\left\{\Big\vert \Vert \Phi\Psi x\Vert_{1}-\mathbb{E}\Vert \Phi\Psi x\Vert_{1} \Big\vert \ge t\right\}\\
        =&\mathbb{P}\left\{\Big\vert\sum_{l=1}^{L}\sum_{i=1}^{m}\Big( \vert(\Phi_{l}\Psi_{l}x)_{i} \vert-\mathbb{E}\vert(\Phi_{l}\Psi_{l}x)_{i} \vert\Big) \Big\vert \ge t\right\},
    \end{align}
    where $(\Phi_{l}\Psi_{l}x)_{i} $ is the $i$-th coordinate of the vector $\Phi_{l}\Psi_{l}x$.
    Note that
\begin{align}
    (\Phi_{l}\Psi_{l}x)_{i}=\sum_{j=1}^{d}(\Phi_{l})_{ij}(\Psi_{l}x)_{j},
\end{align}
where $(\Phi_{l})_{ij}$ is the element in the $i$-th row and $j$-th column of matrix $\Phi_{l}$. Hence, by Lemma \ref{Lem_berstein}, we have for $t\ge 0$
\begin{align}
    \mathbb{P}\big\{\vert\sum_{j=1}^{d}(\Phi_{l})_{ij}(\Psi_{l}x)_{j} \vert\ge t   \big\}\le 2\exp\Big(-c_{1}\min\big(\frac{t^{2}}{\Vert \Psi_{l}x\Vert_{2}^{2}}, \frac{t}{\Vert \Psi_{l}x\Vert_{\infty}}  \big)\Big),
\end{align}
where $c_{1}>0$ is a constant depending only on the subexponential norm $\tau$. A direct integration yields that, for $p\ge 1$
\begin{align}\label{Eq_lemma4.1_2}
    \big(\mathbb{E}\vert(\Phi_{l}\Psi_{l}x)_{i} \vert^{p}\big)^{1/p}\le C_{1}p\Vert\Psi_{l}x\Vert_{2},
\end{align}
where $C_{1}>0$ is a constant depending only on $\tau$. Hence, Lemma \ref{Lem_subexponential_p_moment} yields that $(\Phi_{l}\Psi_{l}x)_{i}$ is a subexponential random variable with subexponential norm no more than $C_{2}\Vert\Psi_{l}x\Vert_{2}$, where $C_{2}>0$ is a constant depending only on $\tau$.

Applying Lemma \ref{Lem_berstein} again, we have for $t\ge 0$
\begin{align}
    &\mathbb{P}\left\{\Big\vert\sum_{l=1}^{L}\sum_{i=1}^{m}\Big( \vert(\Phi_{l}\Psi_{l}x)_{i} \vert-\mathbb{E}\vert(\Phi_{l}\Psi_{l}x)_{i} \vert\Big) \Big\vert \ge t\right\}\\
    \le &2\exp\left(-c_{2}\min\big(\frac{t^{2}}{\sum_{l=1}^{L}\sum_{i=1}^{m}\Vert \Psi_{l}x\Vert_{2}^{2}}, \frac{t}{\max_{l}\Vert\Psi_{l}x\Vert_{2} }  \big)   \right),
\end{align}
where $c>0$ is  a constant depending only on $\tau$. Note that
\begin{align}
    \sum_{l=1}^{L}\Vert \Psi_{l}x\Vert_{2}^{2}=\Vert \Psi x\Vert_{2}^{2}=1, \quad \max_{l}\Vert\Psi_{l}x\Vert_{2}\le 1,
\end{align}
which yields the desired result due to \eqref{Eq_Lemma4.1_1}.
\end{proof}

\subsection{The uniform upper bound}
Define $\Omega:=\{x\in\mathbb{C}^{dL}:  \Vert x\Vert_{\mathcal{I}, 0}\le s, \Vert x\Vert_{2}=1\}$, where $\mathcal{I}=\{\mathcal{I}_{1},\cdots, \mathcal{I}_{G}\}$ is the group partition defined in Theorem \ref{Theo_matrix_satisfies_rip}. We next  show an upper bound for $\sup_{y\in\Omega}\Vert \Phi\Psi y\Vert_{1}$.

\begin{lemma}\label{Lem_uniform_upperbound}
     Consider the matrices $\Phi$ and $\Psi$ as defined in Theorem \ref{Theo_matrix_satisfies_rip}, and let $\Omega$ be the set previously defined. For any $t\ge 0$ and $0<\varepsilon<1$, the following holds
    \begin{align}
    \mathbb{P}\{\sup_{y\in \Omega}\Vert \Phi\Psi y\Vert_{1}\ge\frac{m}{1-\varepsilon}( t+C\sqrt{L} ) \}
    \le \big(\frac{9edL}{\varepsilon^{2}gs}\big)^{gs}\exp\big(-cm\min\big(t^{2}, t \big)  \big),
\end{align}
where $c, C>0$ are constants depending only on the subexponential norm $\tau$, $e$ is the exponential constant, and $g=\max_{i\in [G]}\vert \mathcal{I}_{i} \vert $.
\end{lemma}
\begin{proof}
    By virtue of \eqref{Eq_lemma4.1_2}, we have
    \begin{align}
        \mathbb{E}\Vert \Phi\Psi x\Vert_{1}=\sum_{l=1}^{L}\sum_{i=1}^{m}\mathbb{E}\vert (\Phi_{l}\Psi_{l}x)_{i}\vert\le mC_{1}\sum_{l=1}^{L}\Vert \Psi_{l}x\Vert_{2}\le mC_{1}\sqrt{L},
    \end{align}
    where $C_{1}>0$ is a constant depending only on the subexponential norm $\tau$, and $\Psi_{1}, \cdots, \Psi_{L}$ are submatrices of $\Psi$ as defined in the proof of Lemma \ref{Lem_concentration}. Hence, Lemma \ref{Lem_concentration} yields that
    \begin{align}
        \mathbb{P}\{\Vert \Phi\Psi x\Vert_{1}\ge t+mC_{1}\sqrt{L}  \}\le \exp\big(-c_{1}\min\big(\frac{t^{2}}{m}, t \big)  \big).
    \end{align}

Note that, for any $x\in \Omega$, $\Vert x\Vert_{0}\le gs$. Corollary \ref{Coro_epsilon_net} implies that there exists an $\varepsilon$-net $\mathcal{N}$ of $\Omega$ such that
    \begin{align}
        \vert \mathcal{N}\vert\le \big(\frac{9edL}{\varepsilon^{2}gs}\big)^{gs}.
    \end{align}
Hence, taking the union bound, we have for $t\ge 0$
\begin{align}\label{Eq_Lemma4.2_1}
    \mathbb{P}\{\sup_{x\in \mathcal{N}}\Vert \Phi\Psi x\Vert_{1}\ge t+mC_{1}\sqrt{L}  \}\le& \vert \mathcal{N}\vert\sup_{x\in \mathcal{N}}\mathbb{P}\{\Vert \Phi\Psi x\Vert_{1}\ge t+mC_{1}\sqrt{L}  \}\\
    \le& \big(\frac{9edL}{\varepsilon^{2}gs}\big)^{gs}\exp\big(-c_{1}\min\big(\frac{t^{2}}{m}, t \big)  \big).
\end{align}

Note that for any $y\in \Omega$, there exists $x\in\mathcal{N}$ such that 
\begin{align}
    \Vert x-y\Vert_{2}\le \varepsilon,\quad \Vert x-y\Vert_{\mathcal{I}, 0}\le s.
\end{align}
Hence, we have 
\begin{align}
    \Vert \Phi\Psi y\Vert_{1}=\Vert \Phi\Psi (y-x)+\Phi\Psi x\Vert_{1}\le \frac{\Vert \Phi\Psi (y-x)\Vert_{1}}{\Vert y-x\Vert_{2}}\Vert y-x\Vert_{2}+\Vert\Phi\Psi x\Vert_{1},
\end{align}
which yields that
\begin{align}
    \sup_{y\in \Omega}\Vert \Phi\Psi y\Vert_{1}\le \varepsilon \sup_{y\in \Omega}\Vert \Phi\Psi y\Vert_{1}+\sup_{x\in\mathcal{N}}\Vert \Phi\Psi x\Vert_{1}.
\end{align}
By virtue of \eqref{Eq_Lemma4.2_1}, we have for $0<\varepsilon<1$
\begin{align}
    \mathbb{P}\{\sup_{y\in \Omega}(1-\varepsilon)\Vert \Phi\Psi y\Vert_{1}\ge t+mC_{1}\sqrt{L}  \}
    \le \big(\frac{9edL}{\varepsilon^{2}gs}\big)^{gs}\exp\big(-c_{1}\min\big(\frac{t^{2}}{m}, t \big)  \big),
\end{align}
which concludes the proof.  
\end{proof}

\subsection{The uniform lower bound}
Recall that $\Omega=\{x\in\mathbb{C}^{dL}:  \Vert x\Vert_{\mathcal{I}, 0}\le s, \Vert x\Vert_{2}=1\}$. We next  show a lower bound for $\inf_{y\in\Omega}\Vert \Phi\Psi y\Vert_{1}$.
\begin{lemma}\label{Lem_uniform_lowerbound}
    Consider the matrices $\Phi$ and $\Psi$ as defined in Theorem \ref{Theo_matrix_satisfies_rip}, and let $\Omega$ be the set previously defined. Let $\varepsilon_{0}>0$ be a small enough constant. For any $t\ge 0$ and $0<\varepsilon\le \varepsilon_{0}$,  the following holds
    \begin{align}
&\mathbb{P}\{\inf_{y\in\Omega} \Vert\Phi\Psi y\Vert_{1}\le m(C-2t-\varepsilon C_{0}\sqrt{L})\}\\
\le& 
    2\big(\frac{9edL}{\varepsilon^{2}gs}\big)^{gs}\exp\big(-cm\min(t^{2}, t)  \big),
\end{align}
where $C, C_{0}, c$ are positive constants depending only on the subexponential norm.
\end{lemma}
\begin{proof}
    Recall that 
    \begin{align}
        \mathbb{E}\Vert \Phi\Psi x\Vert_{1}=\sum_{l=1}^{L}\sum_{i=1}^{m}\mathbb{E} \big\vert(\Phi_{l}\Psi_{l} x)_{i}\big\vert.
    \end{align}
   By virtue of \eqref{Eq_lemma4.1_2}, we have for $p\ge 1$
   \begin{align}\label{Eq_Lemma4.3_1}
       \big(\mathbb{E} \big\vert(\Phi_{l}\Psi_{l} x)_{i}\big\vert^{p}\big)^{1/p}\le C_{1} \Vert \Psi_{l}x\Vert_{2}p,
    \end{align}
    where $C_{1}>0$ is an absolute constant. H\"{o}lder's inequality yields that
    \begin{align}
        \mathbb{E}\big\vert (\Phi_{l}\Psi_{l}x)_{i} \big\vert^{2}\le \big( \mathbb{E}\big\vert (\Phi_{l}\Psi_{l}x)_{i} \big\vert\big)^{2/3}\big( \mathbb{E}\big\vert (\Phi_{l}\Psi_{l}x)_{i} \big\vert^{4}\big)^{1/3}.
    \end{align}
    Hence, we have by \eqref{Eq_Lemma4.3_1}
    \begin{align}
        \mathbb{E}\big\vert (\Phi_{l}\Psi_{l}x)_{i} \big\vert\ge \sqrt{\frac{(\mathbb{E}\vert(\Phi_{l}\Psi_{l}x)_{i}) \vert^{2})^{3}}{\mathbb{E}\vert(\Phi_{l}\Psi_{l}x)_{i}) \vert^{4}}}\ge \sqrt{\frac{\Vert \Psi_{l}x\Vert_{2}^{6}}{(4C_{1})^{4}\Vert \Psi_{l}x\Vert_{2}^{4}}}=\frac{\Vert\Psi_{l}x\Vert_{2}}{16C_{1}^{2}},
    \end{align}
    which means that for any $\Vert x\Vert_{2}=1$
    \begin{align}
        \mathbb{E}\Vert \Phi\Psi x\Vert_{1}\ge \frac{m}{16C_{1}^{2}}\sum_{l=1}^{L}\Vert \Psi_{l}x\Vert_{2}\ge \frac{m}{16C_{1}^{2}},
    \end{align}
    where the last inequality is due to $\sum_{l=1}^{L}\Vert \Psi_{l}x\Vert_{2}\ge (\sum_{l=1}^{L}\Vert \Psi_{l}x\Vert_{2}^{2})^{1/2}=1$.
Hence, by Lemma \ref{Lem_concentration}, we have for $t\ge 0$
\begin{align}
    \mathbb{P}\big\{ \Vert\Phi\Psi x\Vert_{1}\le \frac{m}{16C_{1}^{2}}-t  \big\}\le \exp\big(-c_1\min(\frac{t^{2}}{m}, t)  \big).
\end{align}

Corollary \ref{Coro_epsilon_net} yields that, there exists an $\varepsilon$-net $\mathcal{N}$ of $\Omega$ such that
    \begin{align}
        \vert \mathcal{N}\vert\le \big(\frac{9edL}{\varepsilon^{2}gs}\big)^{gs}.
    \end{align}
Taking the union bound, we have
\begin{align}\label{Eq_Lemma4.3_2}
    \mathbb{P}\big\{\inf_{x\in\mathcal{N}} \Vert\Phi\Psi x\Vert_{1}\le \frac{m}{16C_{1}^{2}}-t  \big\}&\le\vert \mathcal{N}\vert \sup_{x\in\mathcal{N}} \mathbb{P}\big\{ \Vert\Phi\Psi x\Vert_{1}\le \frac{m}{16C_{1}^{2}}-t  \big\}\nonumber\\    
    &\le\big(\frac{9edL}{\varepsilon^{2}gs}\big)^{gs}\exp\big(-c_1\min(\frac{t^{2}}{m}, t)  \big).
\end{align}

Note that for any $y\in \Omega$, there exists $x\in\mathcal{N}$ such that 
\begin{align}
    \Vert x-y\Vert_{2}\le \varepsilon,\quad \Vert x-y\Vert_{\mathcal{I}, 0}\le s.
\end{align}
Hence, we have 
\begin{align}
    \Vert \Phi\Psi y\Vert_{1}=\Vert \Phi\Psi (y-x)+\Phi\Psi  x\Vert_{1}\ge \Vert \Phi\Psi 
 x\Vert_{1}-\Vert \Phi\Psi (y-x)\Vert_{1},
\end{align}
which means that
\begin{align}
    \inf_{y\in\Omega}\Vert \Phi\Psi  y\Vert_{1}\ge \inf_{x\in\mathcal{N}}\Vert \Phi\Psi  x\Vert_{1}-\sup_{y\in\Omega}\varepsilon\Vert \Phi\Psi y\Vert_{1}.
\end{align}
Combining with \eqref{Eq_Lemma4.3_2}, we have for $u\ge 0$
\begin{align}
     &\mathbb{P}\{\inf_{y\in\Omega} \Vert\Phi\Psi y\Vert_{1}\le u\}\\
     \le &\mathbb{P}\big\{\inf_{x\in\mathcal{N}} \Vert\Phi\Psi y\Vert_{1}\le \sup_{y\in\Omega}\varepsilon\Vert \Phi\Psi y\Vert_{1}+u  \big\}\\
     \le &\mathbb{P}\big\{\inf_{x\in\mathcal{N}} \Vert\Phi\Psi y\Vert_{1}\le \frac{m\varepsilon}{1-\varepsilon}(t+C_{2}\sqrt{L})+u  \big\}+\mathbb{P}\{\sup_{y\in \Omega}\Vert \Phi\Psi y\Vert_{1}\ge\frac{m}{1-\varepsilon}( t+C_{2}\sqrt{L} ) \},
\end{align}
where $C_{2}$ is the constant appeared in Lemma \ref{Lem_uniform_upperbound}. Let
\begin{align}
    u=\frac{m}{16C_{1}^{2}}-mt-\frac{m\varepsilon}{1-\varepsilon}(t+C_{2}\sqrt{L}).
\end{align}
Using Lemma \ref{Lem_uniform_upperbound} and \eqref{Eq_Lemma4.3_2}, we have for $t\ge 0$
\begin{align}
&\mathbb{P}\{\inf_{y\in\Omega} \Vert\Phi\Psi y\Vert_{1}\le \frac{m}{16C_{1}^{2}}-mt-\frac{m\varepsilon}{1-\varepsilon}(t+C_{2}\sqrt{L})\}\\
\le& 
    2\big(\frac{9edL}{\varepsilon^{2}gs}\big)^{gs}\exp\big(-c_2m\min(t^{2}, t)  \big).
\end{align}
    The desired result follows from adjusting the constants.
\end{proof}

Now, we are prepared to prove Theorem \ref{Theo_matrix_satisfies_rip}.
\begin{proof} [Proof of Theorem \ref{Theo_matrix_satisfies_rip}]
Lemmas \ref{Lem_uniform_upperbound} and \ref{Lem_uniform_lowerbound} yield that, with probability at least
\begin{align}
    1-3\big(\frac{9edL}{\varepsilon^{2}gs}\big)^{gs}\exp\big(-c_1m\min(t^{2}, t)  \big),
\end{align}
    we have for $t\ge 0$ and $0<\varepsilon\le \varepsilon_{0}$ ($\varepsilon_{0}$ is the constant appeared in Lemma \ref{Lem_uniform_lowerbound})
    \begin{align}
        C_{1}-2t-\varepsilon C_{2}\sqrt{L}\le \inf_{y\in\Omega}\frac{1}{m} \Vert\Phi\Psi y\Vert_{1}\le \sup_{y\in\Omega} \frac{1}{m}\Vert\Phi\Psi y\Vert_{1}\le  2t+C_{3}\sqrt{L}. 
    \end{align}

  Let
    \begin{align}
        t=\frac{C_{1}}{8},\quad \varepsilon=\min(\varepsilon_{0}, \frac{C_{1}}{4C_{2}\sqrt{L}}).
    \end{align}
    Hence, with probability at least
    \begin{align}
        1-\big(\frac{C_{4}dL^{2}}{gs}\big)^{gs}e^{-c_{2}m},
    \end{align}
    we have 
    \begin{align}
        C_{4}\le  \inf_{y\in\Omega}\frac{1}{m} \Vert\Phi\Psi y\Vert_{1}\le \sup_{y\in\Omega} \frac{1}{m}\Vert\Phi\Psi y\Vert_{1}\le C_{5}\sqrt{L},
    \end{align}
    where $C_{4}, C_{5}$ are positive constants depending only on the subexopnential norm.

    To ensure that the tail probability tends to $1$, we need to assume that
    \begin{align}
        gs\log\big(\frac{C_{4}dL^{2}}{gs}  \big)\le \frac{1}{2}c_{2}m.
    \end{align}
We complete the proof by adjusting the constants.
\end{proof}
\section{Proof of the main result}
In this section, we shall complete the proof of 
our main result. Before we start the proof, it is recommended that the reader review some of the notations mentioned above.
\begin{proof}[Proof of Theorem \ref{Theo_main}]
    Theorem \ref{Theo_matrix_satisfies_rip} yields that, under the condition 
    \begin{align}
        m\ge C_{1}\max\Big( gLs\log\big(\frac{dL}{gs} \big), \log \eta^{-1}\Big),
    \end{align}
with the probability at least $1-\eta$, we have
\begin{align}
    C_{2}\Vert y\Vert_{2}\le \frac{1}{m}\Vert\Phi\Psi y\Vert_{1}\le C_{3}\sqrt{L}\Vert y\Vert_{2},\quad \forall y\in\Sigma_{\mathcal{I}, c_1 Ls}(\mathbb{C}^{dL}),
\end{align}
where $C_{1}, C_{2}, C_{3}, c_{1}$ are positive constants depending only on the subexponential norm $\tau$, and  it is also required that $c_{2}L>1+C_{3}^{2}L/C_{1}^{2}$. Hence, we conclude the proof by using Theorem \ref{Theo_rip_12}.
\end{proof}

\textbf{Acknowledgment:} This work was supported by the National Key R\&D Program of China (No.2024YFA1013501), the National Natural Science Foundation of China  (No. 12571162), and the Shandong Provincial Natural Science Foundation (No. ZR2024MA082).
\printbibliography

\end{document}